\documentclass[10pt]{article}
\usepackage{latexsym}
\usepackage{fullpage}
\usepackage{amssymb}
\begin{document}
%
% Fonts
%
\font\tenbf=cmbx10  
\font\sevenbf=cmbx7
\font\fivebf=cmbx5
\newfam\bffam
\def\bold{\fam\bffam\tenbf}
\textfont\bffam=\tenbf
\scriptfont\bffam=\sevenbf
\scriptscriptfont\bffam=\fivebf
\font\tengoth=eufm10  
\font\sevengoth=eufm7
\font\fivegoth=eufm5
\newfam\gothfam
\def\goth{\fam\gothfam\tengoth}
\textfont\gothfam=\tengoth
\scriptfont\gothfam=\sevengoth
\scriptscriptfont\gothfam=\fivegoth
\font\bfgraec=cmmib10
%
% Latex-spezifische Abk\"urzungen
%
\def\be{\begin{equation}}
\def\ee{\end{equation}}
\def\bea{\begin{eqnarray}}
\def\eea{\end{eqnarray}}
\def\bdm{\begin{displaymath}}
\def\edm{\end{displaymath}}
%
% einzelne Zeichen
%
\def\a{\alpha}
\def\g{\gamma}
\def\G{\Gamma}
\def\d{\delta}
\def\D{\Delta}
\def\p{\varphi}
\def\t{\tau}
\def\s{\sigma}
\def\S{\Sigma}
\def\k{\kappa}
\def\l{\lambda}
\def\L{\Lambda}
\def\e{\beta}
\def\o{\omega}
\def\O{\Omega}
\def\z{\zeta}
\def\ra{\rangle}
\def\la{\langle}
\def\rb{\rbrack}
\def\lb{\lbrack}
\def\wt{\widetilde}
\def\A{\ifmmode \widehat{\bold A}
        \else $\widehat{\bold A}$ \fi}
\def\nb{\nabla}
\def\part#1{{\partial \over \partial x^{#1}}}
\def\qed{\ifmmode \quad\Box
         \else $\quad\Box$ \fi}
%
% Math. Objekte
%
\def\K#1{\ifmmode {\mathbb{K}}^{#1} 
         \else $ \mathbb{K}^{#1}$ \fi}
\def\R{\mathbb{R}}
\def\C{\mathbb{C}}
\def\Q#1{\ifmmode \mathbb{Q}^{#1}\else $ \mathbb{Q}^{#1}$ \fi}
\def\Z#1{\ifmmode \mathbb{Z}^{#1}\else $ \mathbb{Z}^{#1}$ \fi}
\def\N{\ifmmode \mathbb{N}\else $ \mathbb{N}$ \fi}
\def\H{\mathbb{H}}
\def\RP#1{\ifmmode \R{}P^{#1}
          \else $\R{}P^{#1}$ \fi}
\def\CP#1{\ifmmode \C{}P^{#1}
          \else $\C{}P^{#1}$ \fi}
\def\HP#1{\ifmmode \H{}P^{#1}
          \else $\H{}P^{#1}$ \fi}
\def\KP#1{\ifmmode \K{}P^{#1}
          \else $\K{}P^{#1}$ \fi}
\def\spb{{\bold S}}
\def\spcb{\ifmmode \spb^c
     \else $ \spb^c $ \fi}
%
% Gruppen, Algebren
%
\def\SO{{\bold SO}}
\def\GL{{\bold GL}}
\def\O{{\bold O}}
\def\SU{{\bold SU}}
\def\U{{\bold U}}
\def\PU{{\bold PU}}
\def\Sp{{\bold Sp}}
\def\SL{{\bold SL}}
\def\spin{{\goth spin}}
\def\so{{\goth so}}
\def\su{{\goth su}}
\def\sp{{\goth sp}}
\def\sl{{\goth sl}}
\def\u{{\goth u}}
\def\zg{{\goth g}}
\def\zk{{\goth k}}
\def\zp{{\goth p}}
\def\Spin{{\bold Spin}}
\def\Spinc{\ifmmode {\bold Spin}^c
       \else ${\bold Spin}^c$ \fi}
\def\Spinh{\ifmmode {\bold Spin}^h
       \else ${\bold Spin}^h$ \fi}
\def\Cl{{\bold Cl}}
\def\CCl{\mathbb{C}{\bold l}}
\def\End{{\rm End}}
\def\Hom{{\rm Hom}}
\def\Aut{{\rm Aut}}
\def\Mat{{\rm Mat}}
\def\Bil{{\rm Bil}}
\def\Iso{{\rm Iso}}
\def\Sym{{\rm Sym}}
%
% Abbildungen u.\"a.
%
\def\zff{{\rm II}}
\def\im{{\rm im}}
\def\rg{{\rm rg}}
\def\id{{\rm id}\,}
\def\im{{\rm im}}
\def\ker{{\rm ker}}
\def\coker{{\rm coker}}
\def\det{{\rm det}}
\def\deg{{\rm deg}}
\def\sgn{{\rm sgn}}
\def\tr{{\rm tr}}
\def\dim{{\rm dim}}
\def\inf{\mathop{\rm inf}}
\def\ind{{\rm ind}}
\def\mod{{\rm mod}}
\def\grad{{\rm grad}}
\def\diag{{\rm diag}}
\def\trns{\phantom{}^t}
\def\lange{{\rm L\"ange}}
\def\span{{\rm span}}
\def\Ad{{\rm Ad}}
\def\ad{{\rm ad}}
\def\proj{{\rm proj}}
\def\max{{\rm max}}
\def\Ric{{\rm Ric}}
\def\rank{{\rm rank}\,}
\def\spa{{\rm span}}
\def\und{\qquad \mbox{and} \qquad}
%
% Sonstiges
%
\def\tupl#1#2{(#1_1,\dots,#1_{#2})}
\def\cst{{\rm const}}
\def\us{,\hskip-0.4mm,}
\def\os{``}
\def\buchst#1{\char\expandafter\number#1}
\def\addots{\mathinner{\mkern1mu\raise1pt\vbox{\kern7pt\hbox{.}}
  \mkern2mu\raise4pt\hbox{.}\mkern2mu\raise7pt\hbox{.}\mkern1mu}}
%
%  Spezielle
%
\def\kdirac{{\cal D}}
\def\kdiracq{\bar\kdirac}
\def\haupts{\hbox{\bfgraec\char27}}
\def\HB{{\bold H}}
\def\EB{{\bold E}}
%
% Textformatierung
%
\renewcommand{\theequation}{\thesection.\arabic{equation}}
\def\beweis{{\bf Beweis.}\hskip2mm}
\def\proof{{\bf Proof.}\hskip2mm}
\def\definition{{\bf Definition.}\hskip2mm}
\def\bemerkung{{\bf Bemerkung.}\hskip2mm}
\def\remark{{\bf Remark.}\hskip2mm}
\def\beispiel{{\bf Beispiel.}\hskip2mm}
\def\example{{\bf Example.}\hskip2mm}
\def\leer{\vskip 2.7mm \noindent}
\newtheorem{Satz}{Satz}[section]
\newtheorem{Proposition}{Proposition}[section]
\newtheorem{Theorem}{Theorem}[section]
\newtheorem{Corollary}{Corollary}[section]
\newtheorem{Korollar}{Korollar}[section]
\newtheorem{Lemma}{Lemma}[section]
\newtheorem{Conjecture}{Conjecture}[section]
\def\#{\sharp}
\def\b{\flat}
\def\Cr{{\rm Curv}}
\def\es{\lrcorner\;}
\def\ctr{\lrcorner}
\def\esn{\lrcorner_\circ}
\def\dn{\wedge_\circ}
\def\tfrac#1#2{{\textstyle\frac{#1}{#2}}}
\def\eqref#1{(\ref{#1})}
\title{Quaternionic Killing Spinors}
\author{
 W. Kramer\thanks{Supported by the SFB 256 `Nichtlineare partielle
  Differentialgleichungen'},
 U. Semmelmann\thanks{Supported by the Max--Planck--Institut f\"ur
  Mathematik}, 
 G. Weingart\thanks{Supported by the SFB 256 `Nichtlineare partielle
  Differentialgleichungen'}\\
 {\small Mathematisches Institut der Universit\"at Bonn}\\
 {\small Beringstra\ss{}e 1, 53115 Bonn, Germany}\\
 {\small\texttt{kramer@math.uni-bonn.de,
   uwe@math.uni-bonn.de, gw@math.uni-bonn.de}}}
\maketitle
\begin{abstract}
In \cite{qklast} we proved a lower bound for the spectrum of the Dirac 
operator on quaternionic K\"ahler manifolds. In the present article we
study the limiting case, i.~e.~manifolds where the lower bound is
attained as an eigenvalue. We give an equivalent formulation in terms
of a quaternionic Killing equation and show that the only symmetric
quaternionic K\"ahler manifolds with smallest possible eigenvalue are
the quaternionic projective spaces.
\end{abstract}

\tableofcontents

\section{Introduction}
Let $(M^{4n}, g),\,n\geq 2$ be a compact quaternionic K\"ahler manifold
of positive scalar curvature $\k$. By definition its holonomy group is
then contained in the subgroup $\Sp(n)\cdot\Sp(1)\subset\SO(4n)$. If the
quaternionic dimension $n$ is even or if $M=\HP n$, $M$ is spin and we
proved in \cite{qklast} the following lower bound for the spectrum of
the Dirac operator $D$ on $M$:
$$
\lambda^2 \ge \frac{\k}{4}\,\frac{n+3}{n+2}\,. 
$$
Note that $\k$ is constant on $M$, since any quaternionic K\"ahler
manifold is automatically Einstein. This was first shown by
D.~V.~Alekseevskii in \cite{alex2} and \cite{alex2a} (see also
\cite{ish3}).  The estimate is sharp since the lower bound is the
first eigenvalue of $D^2$ on the quaternionic projective space, as
follows from the computation of the spectrum done in \cite{mil}.

The natural task is then to study the limiting case and find all
manifolds where $\frac{\k}{4}\,\frac{n+3}{n+2}$ is in the spectrum
of $D^2$. In this article we rule out all Wolf spaces besides the
quaternionic projective spaces, thus settling the question for all
compact symmetric quaternionic K\"ahler manifolds. Up to now,
no other examples of compact quaternionic K\"ahler manifolds of
positive scalar curvature are known, and a common conjecture,
proved by C.~LeBrun and S.~Salamon \cite{lebrun} in quaternionic
dimensions $n=2$ and $n=3$, says that there are none.

The principal result shows that the existence of an eigenspinor
with the minimal eigenvalue is equivalent to the existence of a
solution for a suitable quaternionic Killing equation, i.~e.~a
section of a suitable vector bundle which is parallel with respect
to a modified connection. The curvature of this Killing connection
is precisely the hyperk\"ahler or Weyl part of the curvature tensor.
Explicit calculation then shows that no Wolf space besides the
quaternionic projective space allows a parallel section for this
connection.

A peculiar feature of the quaternionic Killing connection is that
unlike its Riemannian or K\"ahlerian counterpart it is not defined
on (a subbundle of) the spinor bundle, but involves a non--spinor
bundle naturally. These ``hidden parameters'' account for the fact
that the dimension of the space of eigenspinors with minimal eigenvalue
on the quaternionic projective space exceeds the dimension of
$\spb_0(\HP n)\oplus\spb_1(\HP n)$. As the geometric significance of
the additional bundle is obscure it seems difficult to describe the
Killing connection in purely geometric terms without using
representation theory of $\Sp(n)\cdot\Sp(1)$.

The authors would like to thank D.~V.~Alekseevskii and C.~B\"ar
for many helpful discussions. Especially, we would like to thank
W.~Ballmann for encouragement and support.

\section{Spin Geometry of Quaternionic K\"ahler Manifolds}

Let $ ( M^{4n}, \, g) $ be a quaternionic K\"ahler manifold,
i.~e.~the Levi--Civit\'a connection on $M$ is already defined on
a $\Sp(n)\cdot\Sp(1)$--reduction $P$ of the  $ \SO(4n)$--bundle
of orthonormal frames. Any representation $V$ of $\Sp(n)\times\Sp(1)$
locally gives a vector bundle ${\bold V}$ associated to (a local
two--fold covering of) $P$. This bundle exists globally iff the
representation factors through $ \Sp(n) \cdot \Sp(1) $.

The representation theory of $\Sp(1)$ and $\Sp(n)$ is governed by
the defining representations $ H : = \H \cong \C^2 $ and
$ E := \H^n \cong \C^{2n} $ respectively. More precisely,
any irreducible $ \Sp(n)\times \Sp(1)$--representation
can be realized as a subspace of $ H^{\otimes p} \otimes E^{\otimes q} $
for some $p$ and $q$; those with $ p + q $ even factor through 
$ \Sp(n) \cdot \Sp(1) $.
Hence, any vector bundle on $ M $ associated to $ P $ can be expressed
in terms of the local bundles $ \HB $ and  $ \EB $. For example, the
complexified tangent bundle is defined by the representation
$ H \otimes E $, i.~e.~
$$ TM^{\C} \; = \; \HB \otimes \EB\, .$$

In this section we will recall some of the definitions and results
given in \cite{qklast}. Besides elementary properties of the 
representations $\L^sE$ and $\Sym^r H$ we also need the explicit
description of the curvature tensor given in \cite{qklast}
as well as the definition of Dirac and twistor operators.

\subsection{Preliminaries on \textbf{Sp}$(n)$--Representations}

Let $H$ and $E$ be the defining complex representations of
$\Sp(1)$ and $\Sp(n)$ with their invariant symplectic forms 
$ \sigma_H \in \L^2 H^\ast$ and $\sigma_E \in \L^2
E^\ast$ and their compatible positive quaternionic structures
$J$, e.~g.~
$$
\begin{array}{rcl}
 J^2  & = & -1 \\
 \s_E(Je_1,Je_2) & = & \overline{\s_E(e_1,e_2)} \\
 \s_E(e,Je) & > & 0 \quad \hbox{for $e\neq 0$\, .}
\end{array}
$$
The symplectic form $\s_E$ defines an isomorphism $\#: E \to E^*,\;
e\mapsto e^\#:=\s_E(e,\,.\,)$ with inverse $\b: E^* \to E $.

Of particular importance is the representation $\Sym^2 E$. Its real
subspace is canonically isomorphic to $\sp(n)$. Thus $\Sym^2 E$ acts
on every complex representation of $\Sp(n)$, e.~g.~its action on
$E$ is given by $(e_1e_2) e := \s_E(e_1, e)e_2 + \s_E(e_2, e)e_1$.
Analogous statements are true for $H$.

Let $ \{ e_i \} $  and $ \{ de_i \} $ with $ de_i ( e_j ) = \delta_{ij} $
be a dual pair of bases for $E$, $E^*$ respectively. In terms of this
bases the symplectic form and its canonical bivector
--- associated by $ \L^2 E \cong \L^2 E^* $ --- can be written as
$$
 \s_E = \tfrac{1}{2} \,\sum de_i \,\wedge \,e_i^\#\,\in \L^2 E^*
 \qquad \qquad
 L_E = \tfrac{1}{2} \,\sum de_i^\b \,\wedge \,e_i\,\in \L^2 E \, .
$$
Wedging with $ L_E $ determines a homomorphism
$ L  : \L^{s-2} E \longrightarrow \L^{s} E$
whereas contracting with $\s_E$ defines its adjoint
$ \L := L^\ast:\, \L^s E \longrightarrow \L^{s-2} E $.
The operators $L$, $\L$ and $H:=\lb \L,L \rb$ with $H|_{\L^sE} = (n-s)\id$
satisfy the commutator algebra of the Lie algebra $\sl_2{\C}$.
Therefore, $\L^sE= \im(L) \oplus \ker(\L)$ splits as 
$\Sp(n)$--representation and $\L^s_\circ E :=\ker(\L)$, the
{\it primitive space}, turns out to be irreducible.
Inductively, the complete decomposition of $\L^sE$
is proved to be:
$$
\L^s E \quad = \quad
\bigoplus^{\lb\frac{s}{2}\rb}_{k=0}  \, \L^{s-2k}_\circ  E, \qquad
0 \le s \le n \, .
$$
The primitive space is stable under contraction with elements
of $E^\ast$ but it is not preserved by the wedge product.
Therefore it is necessary to describe the projection 
$e\dn\omega$ of $e\wedge\omega$ onto $\;\L^s_\circ  E $:

\begin{Lemma} \label{proj}
If $\omega \in   \L^s_\circ E $ then
$e^\#\, \lrcorner \, \omega \; \in  \, \L^{s-1}_\circ E$. Furthermore
$$
e \dn \omega =  e \wedge \omega  \;  - \;
 \tfrac{1}{n-s+1} \, L_E \wedge (e^\#\,\es\,\omega) \, .
$$
\end{Lemma}

Summarizing the properties of contraction and modified exterior
multiplication we have:

\begin{Lemma} \label{kom1}
On $\L^s_\circ E$ modified exterior multiplication and contraction
operators satisfy the following anti--commutator relations
$$ \{ \eta_1\es, \eta_2\es \} =  0 \qquad\qquad
   \{ e_1 \dn , e_2 \dn  \} = 0 \qquad\qquad
   \{ \eta \, \es \, , e \, \dn  \,\} =  \eta ( e )  +
   \tfrac{1}{n-s+1} \eta^\b\, \dn  \,  e^\#\es
$$
for arbitrary $\eta,\eta_i \in E^\ast$ and $e,e_i \in E$. In addition
the following variants of number operators are defined:
$$\begin{array}{lclclcl}
  \sum e_i \, \dn  \, d e_i \, \es & = & s\,\id \, & \qquad &
  \sum d e_i \, \es \, e_i\,\dn & = &\tfrac{(2n - s + 2 )(n-s)}{n-s+1}\,\id.
  \end{array}$$
\end{Lemma}

On $H$, there are similar equations which relate contraction and
symmetric product with $h\in H$. However, it is convenient to modify
contraction. For $\a\in H^*$ we define
$ \alpha \, \esn \, : \Sym^r H \rightarrow \Sym^{r-1} H $ by
$ \alpha \, \esn \, := \tfrac{1}{r} \, \alpha \, \lrcorner \,  $.

\begin{Lemma}\label{kom2}
On $\Sym^rH$ symmetric multiplication and modified contraction operators
satisfy the following commutator relations
$$\begin{array}{lclclcl}
 \lb h_1 \cdot , h_2 \cdot \rb & = & 0 & \qquad\; &
 \lb\alpha \esn  , h \cdot \rb & = & -\frac{1}{r+1} \alpha^\b \cdot h^\#\esn 
  \\\\
 \lb \alpha_1 \esn  , \alpha_2 \esn  \rb & = & 0  & &
  \alpha(h)\id & = & h \cdot \alpha \esn  - \alpha^\b \cdot h^\#\esn
 \end{array}$$
for arbitrary $h,h_i\in H$ and $\alpha,\alpha_i\in H^*$. In addition
the following variants of number operators are defined:
$$\begin{array}{lclclcl}
 \sum h_i \cdot dh_i \esn   & = & \id & &
  \sum dh_i \esn  h_i \cdot & = & \frac{r+2}{r+1} \, \id \, .
\end{array}$$
\end{Lemma}

\subsection{The Curvature Tensor}

For later use we need an explicit description of the curvature
tensor which we take from \cite{qklast}. First we recall the
definition of the following $\End(\HB\otimes\EB)$-valued 
2-forms on $\HB\otimes\EB$:
\begin{eqnarray*}
  R^H_{h_1\otimes e_1,h_2\otimes e_2} &=&
  \s_E(e_1,e_2)(h_1h_2 \otimes \id_E) \nonumber\\
  R^E_{h_1\otimes e_1,h_2\otimes e_2} &=&
  \s_H(h_1,h_2)(\id_H \otimes e_1e_2) \nonumber\\
  R^{hyper}_{h_1\otimes e_1,h_2\otimes e_2} &=& \s_H(h_1,h_2)
  (\id_H \otimes {\goth R}_{e_1,e_2}),
\end{eqnarray*}
where ${\goth R}\in\Sym^4 E^*$ induces the endomorphisms
${\goth R}_{e_1,e_2}: e\mapsto{\goth R}(e_1,e_2,e,.)^\b$ of E.
\leer
\begin{Lemma}\label{curvature}
The curvature tensor of quaternionic K\"ahler manifold $M^{4n}$ is given by
$$
  R = -\frac{\k}{8n(n+2)}(R^H + R^E) + R^{hyper} \, ,
$$
where $\k$ is the scalar curvature of $M$ and the symmetric
4--form ${\goth R}$ is necessarily the symmetrisation: 
$$
  {\goth R}(e_1,e_2,e_3,e_4) = \frac{1}{24 \s_H(h_1,h_2)\s_H(h_3,h_4)}
  \sum_{\t\in\,S_4} \la R_{{h_1}\otimes e_{\t 1},{h_2}\otimes e_{\t 2}}
   {h_3}\otimes e_{\t 3},{h_4}\otimes e_{\t 4}\ra \, ,
$$
which is independent of the choice of the $h_i$ as long as 
$\s_H(h_1,h_2)\s_H(h_3,h_4) \neq 0$.
\end{Lemma}

\subsection{Spinor Bundle and Clifford Multiplication }

The spinor module considered as  $\Sp(n) \times \Sp(1)$--representation
splits into a sum  of $n+1$ irreducible components. Hence, the spinor
bundle of a  $4n$--dimensional  quaternionic K\"ahler manifold
decomposes into a sum of $n+1$ subbundles which can be expressed
using the locally defined bundles $\EB$ and $\HB$.

\begin{Proposition}
{\rm \cite{barsal},\cite{hijmil},\cite{wang}}
The spinor bundle $\spb(M)$ of a quaternionic K\"ahler manifold $M$
decomposes as
$$
\spb(M) \; = \; \bigoplus^n_{r=0} \, \spb_r(M)
  \; := \; \bigoplus^n_{r=0} \, \Sym^r \HB \otimes \L^{n-r}_\circ \EB \, .
$$
The rank of the subbundle $\spb_r(M)$ is given by
$$
\rank (\spb_r(M)) = (r+1)\bigg( {2n \choose n-r } - 
 {2n \choose n-r -2}\bigg) \, .
$$
\end{Proposition}

Note that the covariant derivative on $ \spb(M) $ induced by the
Levi--Civit\'a connection on $(M, g)$ respects the decomposition
given above. The following proposition presents the Clifford
multiplication in terms of the E--H--formalism.

\begin{Proposition}\label{cliff1}{\rm \cite{qklast}}
For any tangent vector $ h \otimes e \in \HB \otimes \EB = TM^{\C}$,
the Clifford multiplication
$\mu( h \otimes e) : \spb(M) \rightarrow \spb(M)  $ is given by:
$$
 \mu(h \otimes e) \quad = \quad \sqrt{2} \;
 (
h \cdot \,  \otimes\, e^\# \, \lrcorner \,
\; + \;
h^\# \, \esn \, \otimes e \, \dn  \,)\, .
$$
In particular, the Clifford multiplication maps the subbundle
$\spb_r(M)$ to the sum $\spb_{r-1}(M) \oplus \spb_{r+1}(M)$.
\end{Proposition}
Thus, Clifford multiplication splits into two components:
\begin{eqnarray*}\label{mu}
 \mu^+_- : \quad  TM \,\otimes\,\spb_r(M)\; \longrightarrow\;\spb_{r+1}(M)
 & \quad\mbox{and}\quad &
 \mu^-_+ : \quad  TM\,\otimes\,\spb_r(M)\; 
 \longrightarrow \; \spb_{r-1}(M)  \, ,
\end{eqnarray*}
where
$ \; \mu^+_-(e\,\otimes\, h\, \otimes\,\psi)
     = \sqrt{2} \, (h \cdot \otimes e^\# \, \es \, ) \, \psi \; $
and
$ \; \mu^-_+ (e\,\otimes\, h\, \otimes\,\psi)
     = \sqrt{2} \,
     (h^\#\,\esn\,\otimes e \, \dn  \,) \, \psi $.
We note that this definition makes sense also for $\spb_r(M)$ replaced
by $ \Sym^r\HB \otimes \L^s_\circ\EB $. In this spirit it is possible
to define two operations similar to Clifford multiplication:
$$
\begin{array}{rccc}
 \mu^+_+: & TM \otimes \Sym^r \HB \otimes \L_\circ^s \EB & \longrightarrow & 
  \Sym^{r+1} \HB \otimes \L_\circ^{s+1} \EB \\ 
 & h \otimes e \otimes \psi & \longmapsto &
  {\sqrt 2}\, (h \cdot \, \otimes \,  e \, \dn\,)\psi
\end{array}
$$
and
$$
\begin{array}{rccc}
 \mu^-_-: & TM \otimes \Sym^r \HB \otimes \L_\circ^s \EB & \longrightarrow & 
  \Sym^{r-1} \HB \otimes \L_\circ^{s-1} \EB \\
 & h \otimes e \otimes \psi &\longmapsto &
  {\sqrt 2}\, (h^\#\,\esn\,\otimes\, e^\#\,\es\,)\psi \, .
\end{array}
$$

Using the number operators of Lemmata \ref{kom1} and \ref{kom2}
it is easy to prove the following useful formulas:

\begin{Lemma}\label{summe}
The following relations are satisfied on
$\Sym^{r}\HB\otimes\L^{s}_\circ\EB$:
$$
\begin{array}{rcl}
 \sum\mu^+_+(X_a)\mu^-_-(X_a) & = & \; 2s\\\\
 \sum\mu^+_-(X_a)\mu^-_+(X_a) & = & - 2\,\tfrac{(2n-s+2)(n-s)}{n-s+1}\\\\
 \sum\mu^-_+(X_a)\mu^+_-(X_a) & = & -2s\,\tfrac{r+2}{r+1}\\\\
 \sum\mu^-_-(X_a)\mu^+_+(X_a) & = & \; 2 \,\tfrac{(2n-s+2)(n-s)}{n-s+1}
   \,\tfrac{r+2}{r+1},
\end{array}
$$
where the sum is over a local orthonormal base $\{X_a\}$ of $TM$.
All other combinations of the partial Clifford multiplications vanish
upon summation over $\{X_a\}$.
\end{Lemma}

\subsection{Dirac and Twistor Operators}

In this section we recall the definition of quaternionic
Dirac and twistor operators. For defining these operators
we have to decompose $ TM \, \otimes \, \spb_r(M) $ into
irreducible components and to project the covariant differential
of a spinor onto the different summands. For $ r\neq 0,n$
we have the following decomposition

\bea \label{deco}
TM \, \otimes \, \spb_r(M)
& \cong &
(\HB \otimes \EB) \, \otimes \, (\Sym^{r}\HB \otimes \L^{n-r}_\circ  \EB) 
\\[1.3ex] \nonumber
& \cong & \spb_{r+1}(M) \;  \oplus \; \spb_{r-1}(M) \; \oplus \;
   ( \, S^+_{r} \; \oplus \; S^-_{r} \;\oplus \; V^+_r 
   \;\oplus\; V^-_r \,)\, .
\eea
Here we used the notation 
$S^{\pm}_{r} =  \Sym^{r \pm 1}\HB \otimes \L^{n-r \pm 1}_\circ  \EB  \; $
and
$ \; V^{\pm}_r = \Sym^{r\pm 1}\HB \otimes K^{n-r}\EB $,
where $ K^{n-r}\EB $ is the summand corresponding to the sum of the highest
weights in the decomposition of $ E \otimes \L^{n-r}_\circ  E$.
In the case $ r = 0 \; $ and $r=n\;$ four of the
above summands vanish and we obtain:
\be\label{rzn}
(\HB \otimes  \EB) \, \otimes \, \L^n_\circ \EB 
\; \cong \; \spb_1(M) \, \oplus \, V^+_0  \qquad \mathrm{and}\qquad
(\HB \otimes  \EB) \, \otimes \, \Sym^n\HB 
\; \cong \; \spb_{n-1}(M) \, \oplus \, S^+_n\, .
\ee

The two components of the Clifford multiplication define natural
projections onto the first two summands appearing in the decomposition
(\ref{deco}). The remaining four summands constitute the kernel
of the Clifford multiplication. The projections onto $ S_r^+$ resp.~$S_r^-$
are given by $\mu^+_+$ resp.~$\mu^-_-$ and  we denote the projections
onto $V^\pm$ by $pr_{V^\pm}$. Applying these projectors to the section
$\nb \psi \in \G(TM \otimes \spb(M))$ we get the two components of the
Dirac operator:
\be \label{defdirac}
D^{+}_{-} \,  := \,
\mu^{+}_{-} \circ \nabla : \spb_r(M)\,\longrightarrow\,\spb_{r +1}(M)\,\qquad 
D^{-}_{+} \,  := \,
\mu^{-}_{+} \circ \nabla : \spb_r(M) \, \longrightarrow\,\spb_{r-1}(M)\, ,
\ee
where $D = D^+_- + D^-_+$ is the Dirac operator, and four twistor
operators:
\be
\begin{array}{rclrcl}
D^{+}_{+} &  := &
\mu^{+}_{+} \circ \nabla : \spb_r(M) \, \longrightarrow \,S_{r}^+ \,\qquad &
D^{-}_{-} &  := &
\mu^{-}_{-} \circ \nabla : \spb_r(M) \, \longrightarrow \,S_{r}^- \, \\\\
T^+ &  := &
pr_{V^+} \circ \nabla : \spb_r(M) \, \longrightarrow \, V^+  \qquad &
T^- &  := &
pr_{V^-} \circ \nabla : \spb_r(M) \, \longrightarrow \, V^- \,. 
\end{array}
\ee
The square of the Dirac operator respects the splitting of
the spinor bundle, i.~e.~$D^2:\spb_r(M)\,\longrightarrow\,\spb_r(M)$.
In particular, we have:
$ D^+_- D^+_- = 0 = D^-_+  D^-_+ $.
\leer

According to the definition of the Dirac and twistor operators
by decomposition (\ref{deco}) it is possible to reconstruct the
covariant differential of a spinor with the help of these operators.
As this is a prerequisite for deriving Killing equations we state
the final formula with the help of right inverses $\iota^{\mp}_{\mp}$
for the partial Clifford multiplications $\mu^{\pm}_{\pm}$:

\begin{Lemma}\label{nabla}
$$
\nabla \phi 
= \iota^-_+(D^+_- \phi) + \iota^+_-(D^-_+ \phi)
  + \iota^-_-(D^+_+ \phi) + \iota^+_+(D^-_- \phi )
  + T^+ \phi + T^- \phi
$$
where the embeddings $\iota^{\mp}_{\mp}$ are defined as the
right inverses of $\mu^{\pm}_{\pm}$. With the help of Lemma
\ref{summe} their explicit form is readily established:
$$
\begin{array}{rcclcrccl}
 \iota^-_+:\!\!&\!\spb_{r+1}(M)\!&\!\longrightarrow\!&\!TM\,\otimes\,\spb_r(M)
 &&
 \iota^-_-:\!\!&\!S^+_r\!&\!\longrightarrow\!&\!TM\,\otimes\,\spb_r(M) \\
 \!&\!\phi\!&\!\longmapsto\!&\!-\tfrac{r+2}{2(n+r+3)(r+1)}
  \sum X_a \otimes \mu^-_+ (X_a)\phi
 &&
 \!&\!\phi\!&\!\longmapsto\!&\!\tfrac{1}{2(n-r+1)}
  \sum X_a \otimes \mu^-_- (X_a)\phi\\\\
 \iota^+_-:\!\!&\!\spb_{r-1}(M)\!&\!\longrightarrow\!&\!TM\,\otimes\,\spb_r(M)
 &&
 \iota^+_+:\!\!&\!S^-_r\!&\!\longrightarrow\!&\!TM\,\otimes\,\spb_r(M) \\
 \!&\!\phi\!&\!\longmapsto\!&\!-\tfrac{r}{2(n-r+1)(r+1)}
  \sum X_a \otimes \mu^+_- (X_a)\phi
 &&
 \!&\!\phi\!&\!\longmapsto\!&\!\tfrac{r(r+2)}{2(n+r+3)(r+1)^2}
  \sum X_a \otimes \mu^+_+ (X_a)\phi
\end{array}
$$
where $\{X_a\}$ is a local orthonormal base of $TM$.
\end{Lemma}

\section{Weitzenb\"ock Formulas}

The central result of \cite{qklast} is a Weitzenb\"ock formula in
matrix form which relates two sets of naturally defined 2nd order
differential operators from the spinor bundle to itself. The idea
is to cope with the abundance of natural 2nd order operators
defined for spin manifolds with special holonomy by replacing
the Lichnerowicz Weitzenb\"ock formula of general holonomy by
the linear space of all Weitzenb\"ock formulas adapted to the
holonomy in question.
In the case of quaternionic K\"ahler manifolds the final matrix
equation is an identity of differential operators defined on
sections of the bundles $\Sym^r\HB\otimes\L_\circ^s\EB$ depending
on the quaternionic dimension $n$:
\be \label{last}
\left(
\begin{array}{c}
  -\nb^\ast\nb\psi \\\\ {\k\over 4}\frac{r(r+2)}{n+2}\psi \\\\
  {\k\over 4}\frac{s(2n-s+2)}{n(n+2)}\psi  \\\\ {\cal C}\psi \\\\
  {\cal L}\psi \\\\  0
\end{array}
\right)
= {\cal W}_H(r)\otimes{\cal W}_E(s)\cdot
\left(
\begin{array}{c}
  -\frac{1}{2}(D^+_+)^\ast D^+_+\psi \\\\ \frac{1}{2}D^+_-D^-_+\psi \\\\
   \frac{1}{2}D^-_+D^+_-\psi \\\\ -\frac{1}{2}(D^-_-)^\ast D^-_-\psi \\\\
  -(T^+)^\ast T^+\psi \\\\ (T^-)^\ast T^-\psi
\end{array}
\right),
\ee
where ${\cal W}_H(r)\otimes{\cal W}_E(s)$ is the Kronecker product
of the two matrices
\be\label{wh}
 {\cal W}_H(r) = \pmatrix{
    1 & -{r\over r+1} \cr\cr
    r & {r(r+2)\over r+1}}
 \qquad\mathrm{and}\qquad
 {\cal W}_E(s) = \pmatrix{
   {1\over s+1}        &-{n-s+2\over(n-s+1)(2n-s+3)}               &1 \cr\cr
  -{s\over s+1}        & {(n-s+2)(2n-s+2)\over(n-s+1)(2n-s+3)}     &1 \cr\cr
  -{(n+1)s\over n(s+1)}&-{(n+1)(n-s)(2n-s+2)\over n(n-s+1)(2n-s+3)}& 
    {n-s\over n}}.
\ee
The proof given in \cite{qklast} for $s=n-r$ goes through without
modification in the general case, only Lemma 4.4 of \cite{qklast}
has to be reformulated. We remark that this formula simplifies in
case $r=0$ or $s=0,n$, because some of the operators involved vanish
by definition. Though this Weitzenb\"ock formula is powerful enough
to prove the eigenvalue estimate for the Dirac operator, it turns
out to be insufficient to derive the quaternionic Killing equations
of the next section.

For that purpose we need additional Weitzenb\"ock formulas
between 2nd order differential operators between different
vector bundles, which are not covered by (\ref{last}). Nevertheless,
the basic idea of \cite{qklast} can be applied to derive these
additional formulas.

We consider for $s\geq 2$ the isotypical $\Sym^r\HB\otimes K^{s-1}
\EB$--component of the second covariant differential $\nabla^2\psi\in
\G(\HB\otimes\EB\otimes\HB\otimes\EB\otimes\Sym^r\HB\otimes\L_\circ^s\EB)$
of a section $\psi$ of $\Sym^r\HB\otimes\L_\circ^s\EB$. As this isotypical
component contains four copies unless $r=0$ the resulting formula will
in general relate two sets of four projectors.

As in \cite{qklast} the problem can be split into two parts dealing
with $\Sp(1)$ and $\Sp(n)$--representations only. Representation theory
of $\Sp(1)$, however, is very simple and no arguments beyond \cite{qklast}
are needed. For this reason we briefly recall that two copies of $\Sym^rH$
are contained in $ H\otimes H\otimes \Sym^rH$, but the projectors
onto these two copies are not unique. A first pair of projectors
is obtained by decomposing $H\otimes H$ into irreducibles, which
then act as endomorphisms on $\Sym^rH$:
$$
\begin{array}{rccl}
  pr_{\C}: \!\!& H \otimes H\otimes\Sym^rH &\longrightarrow&
  \Sym^rH\\[1.4ex]
  & h_1 \otimes h_2 \otimes s & \longmapsto & \s_H(h_1,h_2) s,
  \\\\
  pr_{\Sym^2H}: \!\!& H \otimes H\otimes\Sym^rH &\longrightarrow&
  \Sym^rH\\[1.6ex]
  & h_1 \otimes h_2 \otimes s & \longmapsto & (h_1h_2)(s).
\end{array}
$$
We get a second pair of projectors through the operation of $H$
on $\Sym^rH$ by the $H$-part of Clifford multiplication:
$$
\begin{array}{rccccl}
  pr_{-+}: \!\!& H \otimes H \otimes \Sym^rH & \longrightarrow &
  H \otimes \Sym^{r+1}H &\longrightarrow &  \Sym^rH \\[1.6ex]
  & h_1 \otimes h_2 \otimes s & \longmapsto & h_1 \otimes h_2 \cdot s
    & \longmapsto & h_1^\# \esn(h_2 \cdot s), \\\\
  pr_{+-}: \!\!& H \otimes H \otimes \Sym^rH & \longrightarrow &
  H \otimes \Sym^{r-1}H & \longrightarrow & \Sym^rH \\[1.6ex]
  & h_1 \otimes h_2 \otimes s & \longmapsto & h_1 \otimes h_2^\# \esn s
   & \longmapsto & h_1 \cdot(h_2^\# \esn s).
\end{array}
$$
These two pairs of projectors are related by the matrix ${\cal W}_H(r)$:
\be\label{twisth}
 \pmatrix{ pr_{\C} \cr\cr pr_{\Sym^2H} } =
 \pmatrix{ 1 & -\tfrac{r}{r+1} \cr\cr r & \tfrac{r(r+2)}{r+1} }
 \pmatrix{ pr_{-+} \cr\cr pr_{+-} }.
\ee

Turning now to the second part of the problem dealing with
$\Sp(n)$--representations, we have to look at the isotypical
$K^{s-1}E$--component of $E\otimes E\otimes\L_\circ^2E$
containing two copies. Looking at the decomposition
$E\otimes E\simeq\Sym^2E\oplus\L_\circ^2E\oplus\C{}$
we can write down two projectors immediately:
$$
  \begin{array}{rccccl}
  pr_{K\Sym^2E}: \!\!& E \otimes E \otimes \L_\circ^sE &\longrightarrow&
    \Sym^2E \otimes \L_\circ^sE &\longrightarrow & K^{s-1}E  \\[1.4ex]
  &e_1 \otimes e_2 \otimes \phi & & \longmapsto & &
    \widetilde{pr}_K(e_2\otimes e_1^\#\es\phi+e_1\otimes e_2^\#\es\phi)
  \\\\
  pr_{K\L_\circ^2E}: \!\!& E \otimes E \otimes \L_\circ^sE &\longrightarrow&
    \L_\circ^2E \otimes \L_\circ^sE &\longrightarrow & K^{s-1}E  \\[1.4ex]
  &e_1 \otimes e_2 \otimes \phi & & \longmapsto & &
    \widetilde{pr}_K(e_2\otimes e_1^\#\es\phi-e_1\otimes e_2^\#\es\phi).
\end{array}
$$
On the other hand, we can first project $ E \otimes \L_\circ^sE $ onto
$\L_\circ^{s-1}E $ resp.~$K^{s}E$ and then look what the second $E$-factor
can do:
$$
 \begin{array}{rccccl}
  pr_{K-}: & E \otimes E \otimes \L_\circ^sE & \longrightarrow &
    E \otimes \L_\circ^{s-1}E &\longrightarrow & K^{s-1}E  \\[1.4ex]
  &e_1 \otimes e_2 \otimes \phi & & \longmapsto & &
    \widetilde{pr}_K( e_1\otimes e_2^\#\es\phi), \\\\
  pr_{-K}: & E \otimes E \otimes \L_\circ^sE & \longrightarrow &
    E \otimes K^{s}E &\longrightarrow & K^{s-1}E  \\[1.4ex]
  &e_1 \otimes e_2 \otimes \phi & & \longmapsto & &
    \widetilde{pr}_K\big((\id \otimes e_1^\#\es) 
    \widetilde{pr}_K(e_2\otimes\phi)\big).
 \end{array}
$$
The projector $pr_{-K}$ is not yet in its final form. To simplify
its definition an operator identity on $\L_\circ^{s-1}E$ comes in
handy:
\begin{eqnarray*}
 e_1^\#\es e\dn
 & = &
 -e\dn e_1^\#\es +\s_E(e_1,e)+\tfrac{1}{n-s+2}e_1\dn e^\#\es\\
 & = &
 -\tfrac{(n-s+1)(n-s+3)}{(n-s+2)^2}e\dn e_1^\#\es
 +\tfrac{n-s+1}{n-s+2}\s_E(e_1,e)-\tfrac{1}{n-s+2} e^\#\es e_1\dn\,,
\end{eqnarray*}
which is obtained by applying the anticommutator rules of Lemma
\ref{kom1} twice. We remark that by definition the projector
$\widetilde{pr}_K: E\otimes\L_\circ^sE\to K^sE$,
$$
 \widetilde{pr}_K(e\otimes\phi):= e\otimes\phi
 -\tfrac{1}{s+1}\sum_i e_i\otimes de_i\es e\dn\phi
 -\tfrac{n-s+2}{(2n-s+3)(n-s+1)}\sum_i de^\b_i\otimes e_i\dn e^\#\es\phi,
$$
kills elements of the form $\sum_i e_i\otimes de_i\es\phi$ and
$\sum_i de^\b_i\otimes e_i\dn\phi$. With this in mind the projector
$pr_{-K}$ can be made completely explicit:
\begin{eqnarray*}
 \lefteqn{ pr_{-K}(e_1 \otimes e_2 \otimes \phi) } &&\\
  &=&
   \widetilde{pr}_K\Big(
   e_2 \otimes e_1^\#\es\phi 
   - \tfrac{1}{s+1} \sum_i e_i \otimes e_1^\#\es de_i\es e_2\dn\phi
   - \tfrac{n-s+2}{(2n-s+3)(n-s+1)}
     \sum_i de_i^\b \otimes e_1^\#\es e_i\dn e_2^\#\es\phi \Big)\\
  &=&
   \widetilde{pr}_K(e_2 \otimes e^\#_1\es\phi)\\
  & &
   \! -\tfrac{n-s+2}{(2n-s+3)(n-s+1)} \widetilde{pr}_K
    \Big(\!\sum_i de_i^\b\otimes ( -\tfrac{(n-s+1)(n-s+3)}{(n-s+2)^2}
    e_i\dn e_1^\#\es\!+\!\tfrac{n-s+1}{n-s+2}\s_E(e_1,e_i)\!-\!
    \tfrac{1}{n-s+2} e_i^\#\es e_1\dn )e_2^\#\es\phi\Bigr)\\
  &=&
   \widetilde{pr}_K(e_2 \otimes e^\#_1\es\phi)
   -\tfrac{1}{2n-s+3}\widetilde{pr}_K(e_1 \otimes e^\#_2\es\phi).
\end{eqnarray*}
With this simpler form of $pr_{-K}$ the relations
between the projectors become obvious:
\be\label{twiste}
 \pmatrix{ pr_{K\Sym^2E} \cr\cr pr_{K\L_\circ^2E} } =
 \pmatrix{ \tfrac{2n-s+4}{2n-s+3} & 1 \cr\cr -\tfrac{2n-s+2}{2n-s+3} & 1 }
 \pmatrix{ pr_{K-} \cr\cr pr_{-K} }.
\ee

In a final step the differential operators associated to the projectors
have to be identified. This is simple for the right--hand side projectors,
because by definition the associated operators are products of 1st order
differential operators. To write down the result, we have to define two
new 1st order differential operators, which appear naturally in this way:
$$
  \theta^\pm : \G(\Sym^r\HB \otimes K^s\EB) \to
  \G(\Sym^{r\pm 1}\HB \otimes K^{s-1}\EB)
$$
are the composition of the covariant differential
$\nb:\G(\Sym^r\HB \otimes K^s\EB) \to \G((\HB \otimes
\EB) \otimes\Sym^r\HB \otimes K^s\EB)$ with linear maps
$$
 \begin{array}{rcccl}
  (\HB \otimes \EB) \otimes\Sym^r\HB \otimes K^s\EB &\hookrightarrow &
  (\HB \otimes \EB) \otimes \Sym^r\HB \otimes \EB \otimes \L_\circ^s\EB
  & \to & \Sym^{r\pm 1}\HB \otimes K^{s-1}\EB \\[1mm]
  && (h\otimes e) \otimes s \otimes \tilde e\otimes \phi & \mapsto &
   \left\{\begin{array}{c} h\cdot s\\ h^\#\esn s\end{array}\right\}
   \otimes\widetilde{pr}_K(\tilde e\otimes e^\#\es\phi)\,.
 \end{array}
$$
Then the right--hand projectors define the following operator products:
$$
 \begin{array}{rcl}
  pr_{-+}\otimes pr_{K-}(\nb^2\phi) &=& \frac{1}{\sqrt{2}}T^-D^+_-\phi\\[1mm]
  pr_{+-}\otimes pr_{K-}(\nb^2\phi) &=& \frac{1}{\sqrt{2}}T^+D^-_-\phi\\[1mm]
  pr_{-+}\otimes pr_{-K}(\nb^2\phi) &=& \quad\theta^-T^+\phi \\[1mm]
  pr_{+-}\otimes pr_{-K}(\nb^2\phi) &=& \quad\theta^+T^-\phi \,.
 \end{array}
$$
The left--hand projectors provide two new 2nd order differential operators:
$$
\begin{array}{rcl}
  T_{\cal C}\phi &:=& pr_{\Sym^2H}\otimes pr_{K\Sym^2E}(\nb^2 \phi)\\[1.6ex]
  T_{\cal L}\phi &:=& pr_{\C{}}\otimes pr_{K\L_\circ^2E}(\nb^2 \phi) \,,
\end{array}
$$
which are of no particular importance for the time being, and two linear
operators depending only on curvature: $(pr_{\Sym^2H}\otimes
pr_{K\L_\circ^2E})\circ\nabla^2$ and $(pr_{\C}\otimes
pr_{K\Sym^2E})\circ\nabla^2$.

It is easy to see that $(pr_{\Sym^2H}\otimes pr_{K\L_\circ^2E})\circ\nabla^2$
is the zero operator, because the curvature tensor of a quaternionic
K\"ahler manifold takes values only in the complement $\Sym^2H\otimes\C$
of $\Sym^2H\otimes\L_\circ^2E$ in $\Sym^2H\otimes\L^2E$. We now claim that
the operator $(pr_{\C}\otimes pr_{K\Sym^2E})\circ\nabla^2$
is trivial, too.

Of the three summands of the curvature tensor of a quaternionic K\"ahler
manifold according to Lemma \ref{curvature}, only $R^E$ and $R^{hyper}$
take values in $\C\otimes\Sym^2E$. Tracing down the definitions of
$(pr_{\C}\otimes pr_{K\Sym^2E})\circ\nb^2$ the contribution of
$R^{hyper}$ can be written in terms of ${\goth R}\in\Sym^4E^*$
as the following morphism from $\Sym^rH\otimes\L_\circ^sE$ to
$\Sym^rH\otimes K^{s-1}E$:
$$
 (pr_{\C}\otimes pr_{K\Sym^2E})\circ\nabla^2\,( s\otimes\phi)
 =
 s\otimes\Bigl(\sum_{ij}\big(\widetilde{pr}_K\circ
 (de^\b_i\otimes de_j\es+de^\b_j\otimes de_i\es)\circ{\goth R}_{e_i,e_j}\bigr)
 \,\phi\Bigr)\,.
$$
Expanding $\goth R$ in fourth powers ${1\over 24}\a^4$, $\a\in E^*$ the
contribution from $R^{hyper}$ is seen to vanish because already without
projecting $\a^\b\otimes\a\es\a^\b\dn\a\es\phi=0$ for all $\phi$.
The contribution of $R^E$ to $(pr_{\C}\otimes pr_{K\Sym^2E})\circ\nb^2$
can be written down similarly. The resulting homomorphism is immediately
seen to be $\Sp(n)\cdot\Sp(1)$--equivariant and has to vanish,
because there is no non--trivial way to map $\Sym^rH\otimes\L_\circ^sE$
to $\Sym^rH\otimes K^{s-1}E$ equivariantly.

We conclude this section in tensoring equations (\ref{twisth}) and
(\ref{twiste}) to get the following twistor Weitzenb\"ock formula:
$$
\pmatrix{0 \cr\cr T_{\cal C} \cr\cr T_{\cal L} \cr\cr 0} = {\cal W}_{twist}
\pmatrix{{1\over\sqrt{2}}T^-D^+_- \cr\cr {1\over\sqrt{2}}T^+D^-_-
  \cr\cr \theta^-T^+ \cr\cr \theta^+T^-}
\quad\qquad
{\cal W}_{twist} = 
\pmatrix{
 {2n-s+4\over 2n-s+3} & -{r\over r+1}{2n-s+4\over 2n-s+3}
  & 1 & -{r\over r+1} \cr\cr
 r{2n-s+4\over 2n-s+3} & {r(r+2)\over r+1}{2n-s+4\over 2n-s+3}
  & r & {r(r+2)\over r+1} \cr\cr
 -{2n-s+2\over 2n-s+3} & {r\over r+1}{2n-s+2\over 2n-s+3}
  & 1 & -{r\over r+1} \cr\cr
 -r{2n-s+2\over 2n-s+3} & -{r(r+2)\over r+1}{2n-s+2\over 2n-s+3}
  & r & {r(r+2)\over r+1}}.
$$
\begin{Corollary}\label{twist}
The following operator identity holds on sections of the
bundle $\Sym^r\HB \otimes \L_\circ^s\EB$ with $s\geq 2$:
$$
  {1\over\sqrt{2}}{2n-s+4\over 2n-s+3}
  \Big(T^-D^+_- -{r\over r+1}T^+D^-_- \Big)
  + \Big( \theta^-T^+ -{r\over r+1}\theta^+T^-\Big)
  \quad = \quad 0.
$$
This identity is trivially satisfied for $s=1$, because $T^-D^+_-$,
$T^+D^-_-$, $\theta^-T^+$ and $\theta^+T^-$ all vanish separately.
\end{Corollary}

\section{The Quaternionic Killing Equation}

The matrix Weitzenb\"ock formula (\ref{last}) generates a linear space
of operator identities of 2nd order differential operators from a bundle
$\spb_r(M)$ to itself by multiplying it from the left the an arbitrary
row vector. A particularly important identity in this linear spaces
leads to the following key identity of operator norms for any section
$\psi_r$ of $\spb_r(M)$:
\be\label{min}
 \tfrac{(r+2)(n+r+2)}{n+2}\,\tfrac{\k}{4}\,\|\psi_r\|^2 \; = \; 
  -\,\tfrac{r+1}{n-r+1}\,\|D^+_+\psi_r\|^2\;+\;(r+2)\,\|D^-_+\psi_r\|^2\;+\;
  \tfrac{(r+2)(n+r+2)}{n+r+3}\,\|D^+_-\psi_r\|^2\;-\;2(r+1)\,\|T^+\psi_r\|^2.
\ee
Since $D^2$ respects the splitting of $\spb(M)$ into the subbundles
$\spb_r(M)$ an eigenspinor of $D$ can be assumed to be localized in
$\spb_r(M)\oplus\spb_{r+1}(M)$. If $\psi=\psi_r+\psi_{r+1}$ is such
an eigenspinor with eigenvalue $\l$, then $D^+_-\psi_r=\l\psi_{r+1}$,
$D^-_+\psi_r=0$ and identity (\ref{min}) implies for $\psi_r$:
$$
 \tfrac{\k}{4}\tfrac{n+r+3}{n+2}\|\psi_r\|^2
 =\l^2\|\psi_{r+1}\|^2 - \tfrac{r+1}{r+2}\tfrac{n+r+3}{n+r+2}
  \Bigl( \tfrac{1}{n-r+1}\|D^+_+\psi_r\|^2 + 2\|T^+\psi_r\|^2 \Bigr).
$$
In the same vein identity (\ref{min}) implies for $\psi_{r+1}$:
$$
 \tfrac{\k}{4}\tfrac{n+r+3}{n+2}\|\psi_{r+1}\|^2
 =\l^2\|\psi_r\|^2 - \tfrac{r+2}{r+3}
 \Bigl( \tfrac{1}{n-r}\|D^+_+\psi_{r+1}\|^2 + 2\|T^+\psi_{r+1}\|^2 \Bigr).
$$
These identities prove:

\begin{Theorem}\label{estim}
Let $( M^{4n}, \, g) $ be a compact quaternionic K\"ahler
spin manifold of positive scalar curvature $ \k $ and let
$\psi = \psi_r + \psi_{r+1}\in \G(\spb_r(M)\oplus \spb_{r+1}(M))$
be an eigenspinor for $D$ with eigenvalue $\l$. Then
$$
\l^2 \; \ge  \; \frac{\k}{4} \, \frac{n + r + 3}{n + 2}
$$
with equality if and only if $\psi_r$ and $\psi_{r+1}$ are both
minimal in the sense that $D^+_+\psi_r=0=D^+_+\psi_{r+1}$ and
$T^+\psi_r=0=T^+\psi_{r+1}$.
\end{Theorem}

Thus an eigenspinor $\psi$ of $D$ with smallest possible eigenvalue
$\l$ with $\l^2=\tfrac{\k}{4}\tfrac{n+3}{n+2}$ has to be of the form
$\psi=\psi_0+\psi_1\in\G(\spb_0(M)\oplus\spb_1(M))$. As $\psi_0$
is a section of $\spb_0(M)$ its covariant differential $\nabla\psi$
splits in only two pieces $D^+_-\psi$ and $T^+\psi$ according
to Lemma \ref{deco}. Minimality implies $T^+\psi_0=0$ and Lemma
\ref{nabla} reconstructs $\nabla\psi_0$ from $D^+_-\psi_0=\l\psi_1$:
$$
 \nabla\psi_0 \quad = \quad \l\iota^-_+\psi_1
$$
The covariant differential of $\psi_1$ is more complicated as it
splits into six pieces: $D^\pm_\pm\psi_1$ and $T^\pm\psi_1$.
Minimality implies $D^+_+\psi_1=0$ and $T^+\psi_1=0$, and
as part of the Dirac operator $D^+_-\psi_1=0$. Plugging this
into the Weitzenb\"ock formula (\ref{last}) the second, third
and sixth row read:
$$
 \pmatrix{
   \tfrac{\k}{4}\tfrac{3}{n+2}\psi_1 \cr\cr
   \tfrac{\k}{4}\tfrac{(n+3)(n-1)}{n(n+2)}\psi_1 \cr\cr
   0
 }
 \quad = \quad
 \pmatrix{
  \tfrac{3}{2n}             &-\tfrac{9}{4(n+4)}           &\tfrac{3}{2} \cr\cr
  \tfrac{n-1}{2n}           &-\tfrac{3(n+3)}{4(n+4)}      &-\tfrac{1}{2}\cr\cr
  -\tfrac{3(n-1)(n+1)}{2n^2}&-\tfrac{3(n+3)(n+1)}{4n(n+4)}&\tfrac{3}{2n}
 }
 \pmatrix{
   \tfrac{1}{2}D^+_-D^-_+\psi_1 \cr\cr
  -\tfrac{1}{2}(D^-_-)^*D^-_-\psi_1 \cr\cr
   (T^-)^*T^-\psi_1
 }
$$
With $D^+_-D^-_+\psi_1=\tfrac{\k}{4}\tfrac{n+3}{n+2}\psi_1$ this system of
equations may be solved for $(D^-_-)^*D^-_-\psi_1$ and $(T^-)^*T^-\psi_1$
to find:
$$
 (D^-_-)^*D^-_-\psi_1 = \tfrac{\k}{2}\tfrac{(n+4)(n-1)}{n(n+2)}\psi_1
 \qquad\qquad\mathrm{and}\qquad\qquad
 (T^-)^*T^-\psi_1 = 0\,.
$$
But, since we are on a compact manifold this implies $T^-\psi_1=0$.
Independently, $T^-\psi_1=0$ can be shown by applying Corollary \ref{twist}
to $\psi_0$. With the help of Lemma \ref{nabla} the covariant differential
$\nabla\psi_1$ is reconstructed from $D^-_+\psi_1=\l\psi_0$ and $D^-_-\psi_1$:
$$
 \nabla\psi_1\;=\;\l\,\iota^+_-(\psi_0)\;+\;\iota^+_+(D^-_-\psi_1).
$$
This in turn provides a first version of the quaternionic Killing equations:
$$
\begin{array}{rcl}
\nabla_X \psi_0 & = & -\,\tfrac{\l}{n+3} \,\mu^-_+(X)\,\psi_1 \\\\
\nabla_X \psi_1 & = & -\,\tfrac{\l}{4n}  \,\mu^+_-(X)\,\psi_0
                       \;+\;\tfrac{3}{8(n+4)} \,\mu^+_+(X)\, D^-_-\psi_1.
\end{array}
$$
Unfortunately it is not possible to say much about the section $D^-_-\psi_1$.
The idea to overcome this obstacle is to include this special section and
to consider a quaternionic Killing equation for two spinors and an auxiliary
section of the bundle $\L^{n-2}_\circ\EB$. This yields indeed a useful Killing
equation due to the following proposition:
\leer
\begin{Proposition}
$$
 \nabla_X (D^-_- \psi_1)\;=\;
  -\,\tfrac{\k}{4}\,\tfrac{n+4}{n(n+2)}\,\mu^-_-(X)\,\psi_1 .
$$
\end{Proposition}
\leer
\proof
Since $\psi_-:=D^-_-\psi_1$ is a section of $\L^{n-2}_\circ\EB$ its covariant
differential $\nabla\psi_-$ splits into three pieces: $D^+_\pm\psi_-$ and
$T^+\psi_-$. However, from Corollary \ref{twist} applied to $\psi_1$ we
conclude $T^+\psi_-=T^+D^-_-\psi_1=0$. With $D^+_+=-(D^-_-)^*$ we have
in addition:
$$
 (D^+_+)^* D^+_+ \psi_-
  = D^-_-(D^-_-)^*D^-_-\psi_1
  = \tfrac{\k}{2}\tfrac{(n+4)(n-1)}{n(n+2)}\psi_- .
$$
Now the third row of the Weitzenb\"ock formula (\ref{last}) applied to
$\psi_-$ reads
$$
 \tfrac{\k}{4}\tfrac{(n+4)(n-2)}{n(n+2)}\psi_-
 =
 \tfrac{n-2}{2(n-1)}(D^+_+)^*D^+_+\psi_-
 +\tfrac{2(n+4)}{3(n+5)}D^-_+D^+_-\psi_-
 =
 \tfrac{\k}{4}\tfrac{(n+4)(n-2)}{n(n+2)}\psi_-
 +\tfrac{2(n+4)}{3(n+5)}(D^+_-)^*D^+_-\psi_- .
$$
On the compact manifold $M$ this implies $D^+_-\psi_-=0$. Thus the
covariant differential of $\psi_-$ can be reconstructed from
$D^+_+\psi_-=-(D^-_-)^*D^-_-\psi_1=-
\tfrac{k}{2}\tfrac{(n+4)(n-1)}{n(n+2)}\psi_1$
in the spirit of Lemma \ref{nabla}:
$$
 \nabla\psi_-= -\tfrac{\k}{2}\tfrac{(n+4)(n-1)}{n(n+2)}\iota^-_-(\psi_1) .
$$
The proposition follows.
\qed
\leer
Changing slightly the notation we obtain
\begin{Theorem}
Let $(M^{4n}, g)$ be a compact quaternionic K\"ahler spin manifold of
quaternionic dimension $n$ and with positive scalar curvature $\k$.
Let $ \psi = \psi_0 + \psi_{1} \in \G(\spb_0(M)\oplus \spb_1(M))$
be an eigenspinor for the smallest possible eigenvalue $\l$, i.~e.~
$\l^2 =  \tfrac{\k}{4}\,\tfrac{n+3}{n+2}$. Then $\psi_0, \; \psi_1 $
and $ \psi_- := \tfrac{1}{4\l}\,\tfrac{n+3}{n+4} \,D^-_-\psi_1
\in\G(\L^{n-2}_\circ\EB)$ solve the following quaternionic Killing
equation for the parameter $\l$:
$$
\begin{array}{rccl}
 \nabla_X\psi_0 & = & & \;-\,\tfrac{\l}{n+3}\,\mu^-_+(X)\,\psi_1\\\\
 \nabla_X\psi_1 & = &  -\,\tfrac{\l}{4n} \,\mu^+_-(X)\,\psi_0 &
                      +\,\tfrac{3\l}{2(n+3)}\,\mu^+_+(X)\,\psi_-\\\\
 \nabla_X\psi_- & = &  -\,\tfrac{\l}{4n} \,\mu^-_-(X)\,\psi_1 &
\end{array}
$$
Conversely, if the triple $\psi_0,\,\psi_1,\psi_-$ is a
solution of these equations for any $\l\neq 0$ then
$$
\begin{array}{lclcl}
 D^+_-\psi_0\;=\;\l\,\psi_1 & \qquad\qquad &
 D^-_+\psi_1\;=\;\l\,\psi_0 & \qquad\qquad &
 D^+_+\psi_-\;=\; -\,\l\,\tfrac{n-1}{2n}\,\psi_1 \\[1mm]
 T^+  \psi_0\;=\; 0 &&
 D^-_-\psi_1\;=\; 4\l\,\tfrac{n+4}{n+3}\,\psi_- &&
 D^+_-\psi_-\;=\; 0 \\[1mm]
 && D^+_\pm\psi_1\;=\;0 && T^+\psi_-\;=\;0 \\[1mm]
 && T^\pm\psi_1\;=\;0 &&
\end{array}
$$
In particular, $\psi_0+\psi_1$ is an eigenspinor
for the smallest possible eigenvalue.
\end{Theorem}

We will call solutions of the quaternionic Killing equations
{\it quaternionic Killing spinors}. These are sections 
$\psi = (\psi_0, \, \psi_1,\, \psi_-) $ of the bundle 
$$ 
\spb^{Killing}(M)
 \;:=\;
 \spb_0(M)\,\oplus\,\spb_1(M)\,\oplus\,\L^{n-2}_\circ\EB
 \;\cong\;
 \L^{n}_\circ\EB\,\oplus\,(\HB\otimes\L^{n-1}_\circ\EB)
 \,\oplus\,\L^{n-2}_\circ\EB .
$$
As in the Riemannian or the K\"ahler case it is possible to define a Killing
connection $\nabla^{Killing}$ for which quaternionic Killing spinors are
parallel. On sections of the bundle $\spb^{Killing}(M)$ this connection
is given by 
$ \nabla^{Killing}_X := \nabla_X \; + \; A_X $ where $ A_X $
is the following matrix
$$
 A_X \;=\;\left(\begin{array}{ccc}
 0 & \tfrac{\l}{n+3} \, \mu^-_+(X) & 0 \\\\
 \tfrac{\l}{4n}\,\mu^+_-(X) & 0 & -\,\tfrac{3\l}{2(n+3)}\,\mu^+_+(X) \\\\
 0 & \tfrac{\l}{4n} \, \mu^-_-(X) & 0 
\end{array}\right).
$$
\leer
Hence, quaternionic Killing spinors $\psi = (\psi_0,\,\psi_1,\,\psi_-)$
are annihilated by the curvature of the Killing connection,
i.~e.~$R^{Killing}_{X, Y}\psi = 0$ for any vector fields $X,Y$.
The following proposition is crucial for further investigations
of the Killing equations:
\begin{Proposition}
$$
 R^{Killing}\;=\;R^{hyper}
$$
\end{Proposition}
\leer
\proof
As the partial Clifford multiplications are $\Sp(n)\cdot\Sp(1)$-equivariant,
they are parallel with respect to the Levi--Civit\'a connection, and so
is the endomorphism--valued 1-form $A$. Defining for arbitrary
endomorphism--valued 1-forms $B$ and $C$ the endomorphism--valued
2--form $(B\wedge C)_{X,Y}:=B_X\circ C_Y- B_Y\circ C_X$ the usual
formula relating the curvature of $\nabla$ and $\nabla^{Killing}=\nabla+A$
reads
\be\label{kcurv}
 R^{Killing}
 \; = \; R \; + \; [A,A]
 \; = \; R \; + \;{\l^2\over 4n(n+3)}
 \pmatrix{
  \mu^-_+\wedge\mu^+_- & 0 & {6n\over n+3}\mu^-_+\wedge\mu^+_+ \cr\cr
   0 & \mu^+_-\wedge\mu^-_+ - {3\over 2}\mu^+_+\wedge\mu^-_- & 0 \cr\cr
  {n+3\over 4n}\mu^-_-\wedge\mu^+_- & 0 & -{3\over 2}\mu^-_-\wedge\mu^+_+}.
\ee
The entries of this matrix are endomorphism--valued 2-forms, which
can be simplified further using Lemmas \ref{kom1} and \ref{kom2}.
Considering tangent vectors of the form $X=h_1\otimes e_1$ and
$Y=h_2\otimes e_2$ we get:
\begin{eqnarray*}
 (\mu^-_-\wedge\mu^+_-)_{X,Y}
  & = & 2\s_H(h_1,h_2)(e_1^\#\es e_2^\#\es+e_2^\#\es e_1^\#\es)
  \quad = \quad 0 \\
 (\mu^-_+\wedge\mu^+_+)_{X,Y}
  & = & 2\s_H(h_1,h_2)(e_1\dn e_2\dn + e_2\dn e_1\dn )
  \quad = \quad 0
\end{eqnarray*}
and
\begin{eqnarray*}
 (\mu^-_+\wedge\mu^+_-)_{X,Y}
  & = & 2\s_H(h_1,h_2)(e_1\dn e_2^\#\es + e_2\dn e_1^\#\es)
  \quad = \quad 2 R^E_{X,Y}\\
 (\mu^-_-\wedge\mu^+_+)_{X,Y}
  & = & 2\s_H(h_1,h_2)(e_1^\#\es e_2\dn + e_2^\#\es e_1\dn) \\
  & = & -\tfrac{4}{3}\s_H(h_1,h_2)(e_1\dn e_2^\#\es+e_2\dn e_1^\#\es)
  \quad = \quad -\tfrac{4}{3}R^E_{X,Y}
\end{eqnarray*}
with $e_1^\#\es e_2\dn = -e_2\dn e_1^\#\es+\s_E(e_1,e_2)
+{1\over 3}e_1\dn e_2^\#\es$ on $\L_\circ^{n-2}E$. Using
the same argument the calculation of the last matrix entry
reduces to
\begin{eqnarray*}
 \lefteqn{(\mu^+_-\wedge\mu^-_+-\tfrac{3}{2}\mu^+_+\wedge\mu^-_-)_{X,Y}} &&\\
  &\!\!=\!\!&2\Bigl(h_1\cdot h_2^\#\esn\otimes
     (-e_2\dn e_1^\#\es\!+\!\s_E(e_1,e_2)\!-\!e_1\dn e_2^\#\es)
                   -h_2\cdot h_1^\#\esn\otimes
     (-e_1\dn e_2^\#\es\!+\!\s_E(e_2,e_1)\!-\!e_2\dn e_1^\#\es)\Bigr)\\
  &\!\!=\!\!& 2\s_E(e_1,e_2)
       (h_1\cdot h_2^\#\esn + h_2\cdot h_1^\#\esn)\otimes\id
       +2(h_2\cdot h_1^\#\esn-h_1\cdot h_2^\#\esn)\otimes
       (e_1\dn e_2^\#\esn + e_2\dn e_1^\#\esn) \\
  &\!\!=\!\!& 2(R^H+R^E)_{X,Y}
\end{eqnarray*}
using $h_2\cdot h_1^\#\esn-h_1\cdot h_2^\#\esn=\s_H(h_1,h_2)$ on $H$.
From Lemma \ref{curvature} we know the explicit form of the curvature
tensor $R$. Combining this with the above computations of the matrix
entries in formula (\ref{kcurv}) for $R^{Killing}$ we have
$$
 R^{Killing} = R \; + \; 
 \tfrac{\k}{16n(n+2)}
 \left(
 \begin{array}{cccccc}
 2 R^E & 0 & 0 \\\\
 0 & 2(R^H + R^E) & 0 \\\\
 0 & 0 & 2 R^E
 \end{array}
 \right) = R^{hyper}
$$
because $R^H$ annihilates $\L_\circ^n\EB$ and $\L_\circ^{n-2}\EB$.
\qed
\leer
\begin{Corollary}
 The vector bundle $\spb^{Killing}(\HP n)$ is trivial. Any constant
 section projects to a unique eigenspinor on $\HP n$ with minimal
 eigenvalue.
\end{Corollary}

Let $ \psi  = ( \psi_0, \, \psi_1, \, \psi_-) $ be a quaternionic Killing
spinor. We will now show how to construct eigenfunctions of the Laplace
operator as linear combinations of the length functions of the three
components. Similar constructions were already considered in the Riemannian
and the K\"ahler case. The main tool is the following formula which holds
for any section $\psi$ of an hermitean vector bundle with hermitean
connection $\nabla$:
$$
\Delta \, |\psi|^2 \; = \; 2\;\mathrm{Re}\;(\nabla^* \nabla\,\psi,\,\psi)
\; - \; 2 |\nabla \, \psi|^2 .
$$
From this formula and the quaternionic Killing equations it is easy to derive
$$
\Delta
\left(
\begin{array}{c}
|\psi_0|^2 \\\\
|\psi_1|^2 \\\\
|\psi_-|^2 
\end{array}
\right)
\; = \;
\tfrac{2 \l^2}{n+3} \,
\left(
\begin{array}{ccc}
1 & - \, 1 & 0  \\\\\
- \, \tfrac{n+3}{4n} & 1 & -  \tfrac{6(n+4)}{n+3} \\\\
0 & - \, \tfrac{(n+3)(n-1)}{8n^2} & \tfrac{n+4}{n}
\end{array}
\right)
\,
\left(
\begin{array}{c}
|\psi_0|^2 \\\\
|\psi_1|^2 \\\\
|\psi_-|^2 
\end{array}
\right).
$$
Diagonalizing the above matrix leads to
the definition of the following functions
$$
\begin{array}{rcl}
f_0 & := & \quad
  \tfrac{n+3}{4n}\,|\psi_0|^2\;
 +\;|\psi_1|^2\;
 +\;\tfrac{6n}{n+3}\,|\psi_-|^2\\\\
f_1 & := &
 -\,\tfrac{n+3}{4n}\,|\psi_0|^2\;
 +\;\tfrac{1}{n}\,|\psi_1|^2\;
 +\;\tfrac{2(n+4)}{n+3}\,|\psi_-|^2\\\\
f_2 & := & \quad
 \tfrac{n+3}{4n}\,|\psi_0|^2\;
 -\;\tfrac{n+3}{n}\,|\psi_1|^2\;
 +\;\tfrac{6(n+4)}{n-1}\,|\psi_-|^2.
\end{array}
$$
\begin{Proposition}
The functions $f_0, \, f_1, \, f_2 $ are eigenfunctions of the Laplace
operator and the eigenvalues are the first three eigenvalues of the
Laplace operator on the quaternionic projective space. More precisely,
$ f_0 $ is a constant and
$$
\Delta \, f_1 \; = \;  \tfrac{\k}{2n}\,\tfrac{n+1}{n+2}\,f_1
\quad \mbox {and} \quad
\Delta \, f_2 \; = \; \tfrac{\k}{2n}\,\tfrac{2n+3}{n+2}\,f_2.
$$
\end{Proposition}
\leer
That $f_0$ is constant suggests that $\nabla^{Killing}$ is an hermitean
connection  with respect to a modified scalar product on the bundle
$\spb^{Killing}(M)$ defined for two sections
$ \psi = (\psi_0, \, \psi_1, \, \psi_-) $ and 
$ \phi = (\phi_0, \, \phi_1, \, \phi_-) $ by
$$
\la\psi, \, \phi\ra \;  =
 \;\tfrac{n+3}{4n}\,(\psi_0,\,\phi_0)\;
  + \;(\psi_1,\,\phi_1)\;
  + \;\tfrac{6n}{n+3}\,(\psi_-,\,\phi_-).
$$
\begin{Proposition}
The connection $\nabla^{Killing}$ is hermitean with respect to $\la,\ra$.
\end{Proposition}
\leer
\proof
We have $\nabla^{Killing}_X  = \nabla_X + A_X $. Since $\nabla$ is
an hermitean connection on the different components of $\spb^{Killing}(M)$
we only have to prove that $ A_X $ is a skew--hermitean endomorphism on
$\spb^{Killing}(M)$ with respect to $\la,\ra$ for all real tangent vectors
$X$ or equivalently $\,\mathrm{Re}\;\la A_X \psi,\psi\ra = 0\,$ for all
$\psi$. This is straightforward use of the general formulas
(cf.~\cite{qklast}):
$$
(\mu^-_+ (X) \, \psi_1, \, \psi_0) \; = \;
 - \, (\psi_1, \, \mu^+_-(\overline{X}) \, \psi_0) 
\qquad \mbox{and} \qquad
(\mu^+_+ (X) \, \psi_-, \, \psi_1) \; = \;
(\psi_-, \, \mu^-_-(\overline{X}) \, \psi_1). \quad\qed
$$
\leer
The first non--zero eigenvalue of the Laplace operator on a compact
quaternionic K\"ahler manifold of scalar curvature $\k$ is greater
or equal to $\;\tfrac{\k}{2n}\,\tfrac{n+1}{n+2}\;$ and it is a well
known fact (cf.~\cite{alex3} or \cite{lebrun2}) that equality is attained
if and only if the manifold is isometric to the quaternionic projective
space. Hence, the above construction yields the following proposition:

\begin{Proposition}
Let $(M, \, g) $ be a compact quaternionic K\"ahler spin manifold which
admits a quaternionic Killing spinor such that the associated function
$f_1$ is not identically zero. Then $M$ is isometric to the quaternionic
projective space.
\end{Proposition}

\section{The Killing Curvature on Wolf Spaces}

In this section we determine the eigenvalues of the hyperk\"ahler
part of the curvature tensor considered as a symmetric bilinear form
on $\sp(n)$ for the Wolf spaces, i.~e.~the symmetric compact
quaternionic K\"ahler manifolds. In particular we show that this
symmetric bilinear form is always regular except for the quaternionic
projective space and conclude:

\begin{Lemma}\label{full}
Let $M$ be a Wolf space other than $\HP n$, then the values of
$R^{hyper}$ span $\sp(n)$ at every point $p\in M$ in the following sense
$$ \span\{ R^{hyper}_{X,Y}\,\,\mathrm{with}\,\,X,Y\in T_pM\} = \sp(n) $$
\end{Lemma}

Thus there are no parallel sections of $\spb^{Killing}(M)$, because
if $n>2$ there are no $\sp(n)$-invariant elements in the representation
$\L^n_\circ E\oplus(H\otimes\L^{n-1}_\circ E)\oplus\L^{n-2}_\circ E$.
In case $n=2$ there is up to scalar a unique $\sp(n)$-invariant element
and in consequence the curvature $R^{hyper}$ vanishes on sections of the
trivial line bundle $\L^0_\circ\EB\subset\spb^{Killing}(M)$. Nevertheless
no nontrivial section of this line bundle is parallel with respect to
the Killing connection. This shows that there are no quaternionic
Killing spinors on Wolf spaces besides $\HP n$.

\subsection{The Curvature Endomorphism on Symmetric Spaces}

Let $G/_{\displaystyle K}$ be an symmetric space without euclidean
factors, $\zg=\zk\oplus\zp$ the corresponding decomposition of the
semisimple Lie algebra $\zg$ of $G$ into eigenspaces of the
Cartan--involution. The Killing--form $B$ of $\zg$
is non--degenerate on $\zg$ and assumed either negative or positive
definite on $\zp$. The Riemannian metric on $G/_{\displaystyle K}$
is defined accordingly by $g=\mp B$, the upper sign corresponding
to the negative definite or compact case, the lower to the positive
definite or non--compact case. Despite the choice of metric
the inverse isomorphisms $\#:\zp\to\zp^*$ and $\b:\zp^*\to\zp$
are always taken with respect to the Killing--form.

The decomposition of $\zg=\zk\oplus\zp$ allow to define
two partial Killing--forms $B_\zk$ and $B_\zp$ for $X,Y\in\zg$:
$$
 B_\zk(X,Y)=\tr_\zk(ad_X\circ ad_Y) \qquad\qquad
 B_\zp(X,Y)=\tr_\zp(ad_X\circ ad_Y)
$$
with $B=B_\zk+B_\zp$ by definition. It is well known and easy
to prove that $B_\zk$ and $B_\zp$ are symmetric, vanish on
$\zk\times\zp$ and agree on $\zp\times\zp$, in particular
$B_\zk=B_\zp={1\over 2}B$ on $\zp\times\zp$. However, no
such simple relation is true on $\zk\times\zk$. In fact, if
$$
 \zk=\bigoplus_{i=1}^r \zk_i
$$
is the orthogonal decomposition of $\zk$ into (center and)
simple ideals, by Schur's Lemma there exists constants $l_i\in\R$
such that
$$
 B_\zk\vert_{\zk_i\times\zk_i} = B_{\zk_i}
    = l_i\,B\vert_{\zk_i\times\zk_i}
$$
Let $L$ be the endomorphism of $\zk$ defined by
$L\vert_{\zk_i}=l_i\,\id$, then the relation between
the partial Killing--forms on $\zk\times\zk$ can be
expressed as follows:
$$
  B_\zk(K_1,K_2)  = B( L K_1,K_2 ) \qquad\qquad
  B_\zp(K_1,K_2)  = B((\id-L)K_1,K_2)
$$
Next we define the cobracket $\Delta:\zk\to\L^2\zp$ by
$$
  B(\Delta K,X\wedge Y) := B( K, [X,Y] )
   = B( [K,X],Y ) = B( X^\#\es\Delta K,Y )
$$ 
for all $X,Y\in\zp$, in particular $\Delta$ is a Lie algebra
homomorphism of $\zk$ into $\so\,\zp\cong\L^2\zp$. The cobracket
may be used to define an endomorphism of $\zk$ by $\zk\stackrel
\Delta\longrightarrow \L^2\zp\stackrel{[,]}\longrightarrow\zk$.

\begin{Lemma}
 The endomorphism above equals ${\displaystyle {L-\id\over 2}}$.
\end{Lemma}
To prove this lemma one needs to expand the cobracket in terms of
a dual pair of bases $\{E_i\}$ and $\{dE_i\}$ of $\zp$ and $\zp^*$,
namely $\Delta K=\sum_{i<j}B([K,E_i],E_j)dE_i^\b\wedge dE_j^\b$:
\begin{eqnarray*}
 B(\;[\Delta K_1]\;,K_2\;) & = & B(\Delta K_1, \Delta K_2 ) 
 \quad = \quad {1\over 2}\sum_{ij}B([K_1,E_i],E_j)
               B([K_2,dE_i^\b],dE_j^\b)\\
     & = & -{1\over 2}B([K_2,[K_1,E_i]],dE_i^\b)
 \quad = \quad B({L-\id\over 2}K_1,K_2)
\end{eqnarray*}

The extensions of the metric or the Killing--form $g=\mp B$
to $\L^2\zp$ are positive definite and in fact agree. Thus
the curvature endomorphism $\rho:\L^2\zp\to\L^2\zp$ is
uniquely defined by
$$
 B(\rho(X\wedge Y),Z\wedge W) = B(R_{X,Y}Z,W)
   = -B([X,Y],[Z,W]) = -B(\Delta [X,Y],Z\wedge W)
$$
i.~e.~$\rho=-\Delta\circ [,]$. Together with
$[,]\circ\Delta={L-\id\over 2}$ it follows that
$\rho$ is diagonalizable and may be identified
with ${\id-L\over 2}$ on the image $\Delta\zk$ of $\zk$.
Note that due to $B(R_{X,Y}Z,W)=\pm g(R_{Y,X}Z,W)$ the
definition of $\rho$ is the standard one only in the
compact case. With this in mind the scalar curvature
$\k$ of $G/_{\displaystyle K}$ becomes:
$$\pm\k \quad = \quad 2\,\tr_{\L^2\zp}\rho
        \quad = \quad \tr_\zk(\id-L)$$
Combined with the well known formula for the Ricci curvature
\begin{eqnarray*}
 \Ric(X,Y)
  &     = &      \sum_i B(R_{E_i,X}Y,dE^\b_i)
  \quad = \quad -\sum_i B(dE_i^\b,[Y,[X,E_i]]) \\
  &     = &     -B_\zp(X,Y) \quad = \quad \pm {1\over 2}g(X,Y)
\end{eqnarray*}
implying $\k=\pm{\dim\,\zp\over 2}$ this formula relates the
dimensions of $\zp$ and of the $\zk_i$ with the eigenvalues
of $L$:

\begin{Lemma}\label{trace}
 $$\pm\k = {\dim\,\zp\over 2} = \tr_\zk(\id-L)$$
\end{Lemma}

\subsection{The Eigenvalues of the Curvature Endomorphisms}

\begin{Lemma}
 Let $\zg=\zk\oplus\zp$ with $\zk=\sp(1)\oplus\zk^0$ be the
 Lie algebra decomposition of a Wolf space $G/_{\displaystyle K}$
 of quaternionic dimension $n={1\over 4}\dim\,\zp$. Then
 the Killing--forms $B$ of $\zg$ restricted to $\sp(1)$
 and $B_{\sp(1)}$ of $\sp(1)$ are related by
 $$B_{\sp(1)}={2\over n+2}B$$
 i.~e.~the constant $l_{\sp(1)}$ is always ${2\over n+2}$.
\end{Lemma}
\leer
\proof
From Wolf's construction there exists a base $I,J,K$ of
$\sp(1)$ satisfying $[I,J]=2K$ etc.~such that $ad_I$ is
a complex structure on $\zp$. As $\zk^0$ centralizes $\sp(1)$
the trace of $ad_I\circ ad_I$ on $\zg$ is easily computed:
$$B(I,I) = \hbox{trace on $\sp(1)$ + trace on $\zp$} = -8-4n.$$
\noindent On the other hand $B_{\sp(1)}(I,I)=-8={2\over n+2}B(I,I)$. 
\qed
\leer
This result together with Lemma \ref{trace} will be used in the
sequel to calculate the curvature endomorphisms $\rho$ of all
Wolf--spaces:
\leer
{
 \parindent0cm
 \everypar={\hangindent47mm\hangafter0}
 \def\item[#1]{\par\hskip-47mm\smash{\hbox to 0pt{#1:\hss}}\hskip47mm}
 \item[$\Sp(n+1)/_{\displaystyle\Sp(1)\times\Sp(n)}$]
  On the quaternionic projective spaces $\zk$ decomposes into
  $\sp(1)\oplus\sp(n)$. As $\sp(n)$ is the compact form
  of $\sp\,\C^{2n}$ the Killing--form is $(2n+2)\tr_{\C^{2n}}$.
  Now $\sp(n)$ is embedded in $\sp(n+1)$ via the embedding of
  the defining representations $\C^{2n}\hookrightarrow\C^{2n+2}$.
  In consequence, the Killing--forms are related by the constant
  $l_{\sp(n)}={2n+2\over 2n+4}={n+1\over n+2}$. Thus the curvature
  endomorphism of the quaternionic projective space reads:
  $$\rho = {n\over 2(n+2)}\proj_{\sp(1)}
           +{1\over 2(n+2)}\proj_{\sp(n)}.$$
 \item[$\SU(n+2)/_{\displaystyle{\bold S}(\U(2)\times\U(n))}$]
  On the complex Grassmannian of 2-planes $\zk$ decomposes
  into $\sp(1)\oplus\su(n)\oplus\R$. The constant $l_\R$ is
  certainly $0$. As $\su(n)$ is the compact form of $\sl(n)$
  the Killing--form is $2n\;\tr_{\C^n}$ and one concludes that
  $l_{\su(n)}={2n\over 2n+4}={n\over n+2}$. The curvature
  endomorphism of the complex Grassmannian is thus:
  $$\rho = {n\over 2(n+2)}\proj_{\sp(1)}
           +{1\over n+2}\proj_{\su(n)}
           +{1\over 2}\proj_\R.$$
 \item[$\SO(n+4)/_{\displaystyle{\bold S}(\O(4)\times\O(n))}$]
  On the real Grassmannian of 4-planes $\zk$ decomposes into
  $\sp(1)\oplus\widetilde{\sp}(1)\oplus\so(n)$. The Killing--form
  of $\so(n)$ is $(n-2)\tr_{\C^n}$. Arguing in the same way as
  before one finds $l_{\so(n)}={n-2\over n+2}$. The last constant
  is then calculated with the help of Lemma \ref{trace} and agrees
  with $l_{\sp(1)}=l_{\widetilde{\sp}(1)}={2\over n+2}$, because either
  of the two subalgebras could be used to define the quaternionic
  structure. The curvature endomorphism is:
  $$\rho = {n\over 2(n+2)}\proj_{\sp(1)}
           +{n\over 2(n+2)}\proj_{\widetilde{\sp}(1)}
           +{2\over n+2}\proj_{\so(n)}.$$
 \item[${\bold G}_2/_{\displaystyle\SO(4)}$]
  In this case $\dim\,\zp=8$ or $n=2$ and
  $\zk=\sp(1)\oplus\widetilde{\sp}(1)$. However
  $\sp(1)$ and $\widetilde{\sp}(1)$ are definitely
  distinct, $\widetilde{\sp}(1)$ does not define a proper
  quaternionic structure. In fact their constants differ,
  because Lemma \ref{trace} implies
  $l_{\widetilde{\sp}(1)}={1\over 6}$.
  $$\rho = {1\over 4}\proj_{\sp(1)}
           +{5\over 12}\proj_{\widetilde{\sp}(1)}.$$
 \item[${\bold F}_4/_{\displaystyle\Sp(1)\Sp(3)}$]
  This space is not spin, however Lemma \ref{full} may be
  of independent interest. With $\dim\,\zp=28$ the quaternionic
  dimension is $n=7$ and $\zk=\sp(1)\oplus\sp(3)$. Using
  Lemma \ref{trace} one finds $l_{\sp(3)}={4\over 9}$ and
  $$\rho = {7\over 18}\proj_{\sp(1)}
           +{5\over 18}\proj_{\sp(3)}.$$
 \item[${\bold E}_6/_{\displaystyle\Sp(1)\SU(6)}$]
  With $\dim\,\zp=40$ the quaternionic dimension is $n=10$,
  $\zk=\sp(1)\oplus\su(6)$, $l_{\su(6)}={1\over 2}$ and
  $$\rho = {5\over 12}\proj_{\sp(1)}
           +{1\over 4}\proj_{\su(6)}.$$
 \item[${\bold E}_7/_{\displaystyle\Sp(1)\Spin(12)}$]
  With $\dim\,\zp=64$ the quaternionic dimension is $n=16$,
  $\zk=\sp(1)\oplus\so(12)$, $l_{\so(12)}={5\over 9}$ and
  $$\rho = {4\over 9}\proj_{\sp(1)}
           +{2\over 9}\proj_{\so(12)}.$$
 \item[${\bold E}_8/_{\displaystyle\Sp(1){\bold E}_7}$]
  With $\dim\,\zp=112$ the quaternionic dimension is $n=28$
  and $\zk=\sp(1)\oplus{\goth e}_7$, where $\dim\,{\goth e}_7=133$.
  Applying Lemma \ref{trace} one finds $l_{{\goth e}_7}={3\over 5}$
  and thus
  $$\rho = {7\over 15}\proj_{\sp(1)}
           +{1\over 5}\proj_{{\goth e}_7}.$$
}

To present the argument leading to Theorem \ref{full} we
recall that for any Wolf space we have a $K$-invariant subalgebra
$\sp(n)\subset\so\,\zp$ defined as the centralizer of the
subalgebra $\sp(1)\subset\so\,\zp$ defining the quaternionic
structure and by definition $\zk^0\subset\sp(n)$. In particular,
the curvature endomorphism $\rho^{\HP n}$ of the quaternionic
projective space is well defined on any Wolf space.
Bridging the gap between Lie algebra-- and E--H--formalism, it
is even possible to identify the ``curvature'' endomorphisms
corresponding to $R^H$ and $R^E$ of Lemma \ref{curvature} with:
$$
 \rho^H = -2n\,\proj_{\sp(1)}
 \qquad\qquad\mathrm{and}\qquad\qquad
 \rho^E = -2\,\proj_{\sp(n)}
$$
The scalar curvature of a compact Wolf space is $\k=2n$
and Lemma \ref{curvature} implies for the curvature tensor:
$$
 \rho 
 = -{\k\over 8n(n+2)}(\rho^H+\rho^E) + \rho^{hyper}
 = {n\over 2(n+2)}\proj_{\sp(1)}+{1\over 2(n+2)}\proj_{\sp(n)} +\rho^{hyper}.
$$
We conclude that the kernel of the hyperk\"ahler part $\rho^{hyper}$
of the curvature endomorphism $\rho$ restricted to $\sp(n)\subset\so\,\zp$
is the eigenspace of $\rho$ with eigenvalue $\tfrac{1}{2(n+2)}$ in $\sp(n)$.
Looking down the list of curvature endomorphisms above one verifies that
this eigenvalue does only occur for the curvature endomorphism of the
quaternionic projective space itself.
\qed

\end{document}